\documentclass[12pt,twoside]{article}
\usepackage[english]{babel}
\usepackage[latin1]{inputenc}
\usepackage{amsmath}
\usepackage{amssymb,amsfonts}
\usepackage{graphicx}
\usepackage{times,amssymb,amscd}

\newcommand{\bC}{\mathbf{C}}
\newcommand{\bE}{\mathbf{E}}
\newcommand{\bG}{\mathbf{G}}
\newcommand{\bH}{\mathbf{H}}

\newcommand{\bL}{\mathbf{L}}

\newcommand{\bR}{\mathbf{R}}
\newcommand{\bS}{\mathbf{S}}

\newcommand{\bs}{\mathbf{s}}
\newcommand{\ba}{\mathbf{a}}

\newcommand{\bT}{\mathbf{T}}

\newcommand{\bb}{\mathbf{b}}
\newcommand{\bt}{\mathbf{t}}

\newcommand{\cP}{\mathcal{P}}
\newcommand{\cS}{\mathcal{S}}
\newcommand{\cT}{\mathcal{T}}
\newcommand{\cC}{\mathcal{C}}
\newcommand{\cH}{\mathcal{H}}
\newcommand{\cM}{\mathcal{M}}

\newcommand{\EUC}{\bE^3}

\newcommand{\HYP}{\bH^3}
\newcommand{\SXR}{\bS^2\!\times\!\bR}
\newcommand{\HXR}{\bH^2\!\times\!\bR}
\newcommand{\SLR}{\widetilde{\bS\bL_2\bR}}
\newcommand{\NIL}{\mathbf{Nil}}
\newcommand{\SOL}{\mathbf{Sol}}

\newtheorem{Definition}{Definition}[section]
\newtheorem{Remark}{Remark}[section]
\newtheorem{Corollary}{Corollary}[section]

\begin{document}
\pagestyle{myheadings}
\markboth{\centerline{Jen\H o Szirmai}}
{Fibre-like cylinders, their packings...}
\title
{Fibre-like cylinders, their packings and coverings in $\SLR$ space
\footnote{Mathematics Subject Classification 2010: 52C17, 52C22, 52B15, 53A35, 51M20. \newline
Key words and phrases: Thurston geometries, $\SLR$ geometry, cylinder, cylinder packing and covering and their densities, regular
prism tiling \newline
}}

\author{Jen\H o Szirmai \\
\normalsize Department of Algebra and Geometry, Institute of Mathematics,\\
\normalsize Budapest University of Technology and Economics, \\
\normalsize M\"uegyetem rkp. 3., H-1111 Budapest, Hungary \\
\normalsize szirmai@math.bme.hu \\
\date{\normalsize{\today}}}


\maketitle
\begin{abstract}

In this paper we define the notion of infinite or bounded fibre-like geodesic cylinder in $\SLR$ space, develop a method to determine its volume and total surface area.
We prove that the common part of the above congruent fibre-like cylinders with the base plane are Euclidean circles and determine their radii.

Using the former classified infinite or bounded congruent regular prism tilings with generating groups $\mathbf{pq2_1}$ 
we introduce the notions of cylinder packings, coverings and their densities. 
Moreover, we determine the densest packing, the thinnest covering cylinder arrangements in $\SLR$ space, their densities, their connections 
with the extremal hyperbolic circle arrangements and with the extremal fibre-like cylinder arrangements in $\HXR$ space

In our work we use the projective model of $\SLR$ introduced by E. {Moln\'ar} in \cite{M97}.
\end{abstract}

\newtheorem{theorem}{Theorem}[section]
\newtheorem{corollary}[theorem]{Corollary}
\newtheorem{conjecture}[theorem]{Conjecture}
\newtheorem{lemma}[theorem]{Lemma}
\newenvironment{example}{Example}



\section{Introduction}
In \cite{Sz13-1} we defined and described the {\it regular infinite or bounded} $p$-gonal prism tilings in $\SLR$ space.
We proved that there exist infinitely many regular infinite $p$-gonal face-to-face prism tilings $\cT^i_p(q)$ and
infinitely many regular bounded $p$-gonal non-face-to-face prism tilings $\cT_p(q)$ for integer parameters $p,q;~3 \le p$,
$ \frac{2p}{p-2} < q$ with generating groups $\mathbf{pq2_1}$.

In \cite{MSz14} we considered the problem of geodesic ball packings related to tilings and their symmetry groups $\mathbf{pq2_1}$.
Moreover, we computed the volumes of prisms and defined the
notion of geodesic ball packing and its density. In \cite{MSz14} we developed a
procedure to determine the densities of the densest geodesic ball
packings for the tilings considered, more precisely, for
their generating groups $\mathbf{pq2_1}$ (for integer rotational
parameters $p,q$; $3\le p,~\frac{2p}{p-2} <q$). We looked for those
parameters $p$ and $q$ above, where the packing density as largest as possible.
Currently our record in these cases is $0.5674$ for $(p, q) = (8,
10)$. 

In \cite{Sz14} we studied the structure of the regular infinite or bounded $p$-gonal prism tilings and we proved that the side curves of their base figures
are arcs of Euclidean circles for each parameter. Furthermore, we examined the non-periodic geodesic
ball packings of congruent regular non-periodic prism tilings derived from the regular infinite $p$-gonal
face-to-face prism tilings $\cT^i_p(q)$ in $\SLR$ geometry. We determined
the densities of the above non-periodic optimal geodesic ball packings and applied this algorithm to them.
We search for values of parameters $p$ and $q$ that provided the largest packing density.
In this paper we obtain greater density $0.626606\dots$ for $(p, q) = (29,3)$
than the maximum density of the corresponding periodic geodesic ball packings under the groups $\mathbf{pq2_1}$.

In \cite{MSzV13} and \cite{MSz14} we considered tilings $\cT(p,(q,k),(o,\ell))$ for suitable integer positive parameters $p,q,k,o,\ell$.
Every tiling $\cT$ is generated by discrete isometry group $\mathbf{pq}_{k}\mathbf{o}_{\ell}$ for $k=1, o=2$, $\ell=1$.
That means this group is generated by a $p-$rotation $\mathbf{p}$ about the central fibre, then by $\mathbf{q}_k$ screw with $q-$rotation and $\frac{k}{q}$
fibre translation, then by an $\mathbf{o}_\ell$ screw with $o-$rotation and $\frac{\ell}{o}$ translation, just by Euclidean analogy but exact projective computations.
We computed the maximal density of the geodesic and translation ball packings induced by the
$\mathbf{pq}_{k}\mathbf{o}_{\ell}$ group action for any parameters. The highest density in these cases is $\approx 0.787758$ for geodesic packings and 
$\approx 0.84170$ for translation ball packings.

In \cite{CsSz16} we studied the interior angle sums of translation and geodesic triangles in $\SLR$ geometry.
We proved that the angle sum $\sum_{i=1}^3(\alpha_i) \ge \pi$ for translation triangles and for geodesic triangles the angle sum can 
be larger, equal or less than $\pi$. 

In this paper we define the notion of infinite or bounded fibre-like geodesic cylinder in $\SLR$ space, develop a method to determine its volume and total surface area 
(see Lemma 4.1-5 and Corollary 4.1).
We prove that the common part of the above congruent fibre-like cylinders with the base plane are Euclidean circles and determine their radii.

Using the former classified infinite or bounded congruent regular prism tilings with generating groups $\mathbf{pq2_1}$ 
we introduce the notions of cylinder packings, coverings and their densities (see Definition 4.3-4). 

Moreover, we determine the densest packing, the thinnest covering cylinder arrangements in $\SLR$ space, their densities and their connections 
with the extremal hyperbolic circle arrangements and with the extremal fibre-like cylinder arrangements in $\HXR$ space (see Table 2-3, Theorem 5.1-2, Corollary 5.1). 
\section{Preliminaries and background}
The real $ 2\times 2$ matrices $\begin{pmatrix}
         d&b \\
         c&a \\
         \end{pmatrix}$ with unit determinant $ad-bc=1$
constitute a Lie transformation group by the usual product operation, taken to act on row matrices as on point coordinates on the right as follows
\begin{equation}
\begin{gathered}
(z^0,z^1)\begin{pmatrix}
         d&b \\
         c&a \\
         \end{pmatrix}=(z^0d+z^1c, z^0 b+z^1a)=(w^0,w^1)\\
\mathrm{with} \ w=\frac{w^1}{w^0}=\frac{b+\frac{z^1}{z^0}a}{d+\frac{z^1}{z^0}c}=\frac{b+za}{d+zc}, \tag{2.1}
\end{gathered}
\end{equation}
as action on the complex projective line $\bC^{\infty}$ (see \cite{M97}, \cite{MSz}).
This group is a $3$-dimensional manifold, because of its $3$ independent real coordinates and with its usual neighbourhood topology (\cite{R}, \cite{S}, \cite{T}).
In order to model the above structure in the projective sphere $\cP \cS^3$ and in the projective space $\cP^3$ (see \cite{M97}),
we introduce the new projective coordinates $(x^0,x^1,x^2,x^3)$ where
$a:=x^0+x^3, \ b:=x^1+x^2, \ c:=-x^1+x^2, \ d:=x^0-x^3$
with the positive, then the non-zero multiplicative equivalence as projective freedom in $\cP \cS^3$ and in $\cP^3$, respectively.
Then it follows that $0>bc-ad=-x^0x^0-x^1x^1+x^2x^2+x^3x^3$
describes the interior of the above one-sheeted hyperboloid solid $\cH$ in the usual Euclidean coordinate simplex with the origin
$E_0(1;0;0;0)$ and the ideal points of the axes $E_1^\infty(0;1;0;0)$, $E_2^\infty(0;0;1;0)$, $E_3^\infty(0;0;0;1)$.
We consider the collineation group ${\bf G}_*$ that acts on the projective sphere $\cS\cP^3$  and preserves a polarity i.e. a scalar product of signature
$(- - + +)$, this group leaves the one sheeted hyperboloid solid $\cH$ invariant.
We have to choose an appropriate subgroup $\mathbf{G}$ of $\mathbf{G}_*$ as isometry group, then the universal covering group and space
$\widetilde{\cH}$ of $\cH$ will be the hyperboloid model of $\SLR$ \cite{M97}.

The specific isometries $\bS(\phi)$ $(\phi \in \bR )$ constitute a one parameter group given by the matrices:
\begin{equation}
\begin{gathered} \bS(\phi):~(s_i^j(\phi))=
\begin{pmatrix}
\cos{\phi}&\sin{\phi}&0&0 \\
-\sin{\phi}&\cos{\phi}&0&0 \\
0&0&\cos{\phi}&-\sin{\phi} \\
0&0&\sin{\phi}&\cos{\phi}
\end{pmatrix}
\end{gathered} \tag{2.2}
\end{equation}
The elements of $\bS(\phi)$ are the so-called {\it fibre translations}. We obtain a unique fibre line to each $X(x^0;x^1;x^2;x^3) \in \widetilde{\cH}$
as the orbit under the right action of $\bS(\phi)$ on $X$. The coordinates of points lying on the fibre line through $X$ can be expressed
as the images of $X$ by $\bS(\phi)$:
\begin{equation}
\begin{gathered}
(x^0;x^1;x^2;x^3) \stackrel{\bS(\phi)}{\longrightarrow} {(x^0 \cos{\phi}-x^1 \sin{\phi}; x^0 \sin{\phi} + x^1 \cos{\phi};} \\ {x^2 \cos{\phi} + x^3 \sin{\phi};-x^2 \sin{\phi}+
x^3 \cos{\phi})}.
\end{gathered} \tag{2.3}
\end{equation}
The points of a fibre line through
$X$ by usual inhomogeneous Euclidean coordinates $x=\frac{x^1}{x^0}$, $y=\frac{x^2}{x^0}$, $z=\frac{x^3}{x^0}$, $x^0\ne 0$ are given by
\begin{equation}
\begin{gathered}
(1;x;y;z) \stackrel{\bS(\phi)}{\longrightarrow} {\Big( 1; \frac{x+\tan{\phi}}{1-x \tan{\phi}}; \frac{y+z \tan{\phi}}{1-x \tan{\phi}};
\frac{z - y \tan{\phi}}{1-x \tan{\phi}}\Big)}
\end{gathered} \tag{2.4}
\end{equation}
for the projective space $\cP^3$, where ideal points (at infinity) conventionally occur.

In (2.3) and (2.4) we can see the $2\pi$ periodicity of $\phi$, moreover the (logical) extension to $\phi \in \bR$, as real parameter, to have
the universal covers $\widetilde{\cH}$ and $\SLR$, respectively, through the projective sphere $\cP\cS^3$. The elements of the isometry group of
$\mathbf{SL_2R}$ (and so by the above extension the isometries of $\SLR$) can be described by the matrix $(a_i^j)$ (see \cite{M97} and \cite{MSz})
Moreover, we have the projective proportionality, of course.
We define the {\it translation group} $\bG_T$, as a subgroup of the isometry group of $\mathbf{SL_2R}$,
the isometries acting transitively on the points of ${\cH}$ and by the above extension on the points of $\SLR$ and $\widetilde{\cH}$.
$\bG_T$ maps the origin $E_0(1;0;0;0)$ onto $X(x^0;x^1;x^2;x^3)$. These isometries and their inverses (up to a positive determinant factor)
are given by the following matrices:
\begin{equation}
\begin{gathered} \bT:~(t_i^j)=
\begin{pmatrix}
x^0&x^1&x^2&x^3 \\
-x^1&x^0&x^3&-x^2 \\
x^2&x^3&x^0&x^1 \\
x^3&-x^2&-x^1&x^0
\end{pmatrix}.
\end{gathered} \tag{2.5}
\end{equation}
The rotation about the fibre line through the origin $E_0(1;0;0;0)$ by angle $\omega$ $(-\pi<\omega\le \pi)$ can be expressed by the following matrix
(see \cite{M97})
\begin{equation}
\begin{gathered} \bR_{E_O}(\omega):~(r_i^j(E_0,\omega))=
\begin{pmatrix}
1&0&0&0 \\
0&1&0&0 \\
0&0&\cos{\omega}&\sin{\omega} \\
0&0&-\sin{\omega}&\cos{\omega}
\end{pmatrix},
\end{gathered} \tag{2.6}
\end{equation}
and the rotation $\bR_X(\omega)$ about the fibre line through $X(x^0;x^1;x^2;x^3)$ by angle $\omega$ can be derived by formulas (2.5) and (2.6):
\begin{equation}
\bR_X(\omega)=\bT^{-1} \bR_{E_O} (\omega) \bT:~(r_i^j(X,\omega)).
\tag{2.7}
\end{equation}
Horizontal intersection of the hyperboloid solid $\cH$ with the plane $E_0 E_2^\infty E_3^\infty$ provides the
{\it hyperbolic $\mathbf{H}^2$ base plane} of the model $\widetilde{\cH}=\SLR$.
The fibre through $X$ intersects the base plane $z^1=x=0$ in the foot point
\begin{equation}
\begin{gathered}
Z(z^0=x^0 x^0+x^1x^1; z^1=0; z^2=x^0x^2-x^1x^3;z^3=x^0x^3+x^1x^2).
\end{gathered} \tag{2.8}
\end{equation}
After \cite{M97}, we introduce the so-called hyperboloid parametrization as follows
\begin{equation}
\begin{gathered}
x^0=\cosh{r} \cos{\phi}, ~ ~
x^1=\cosh{r} \sin{\phi}, \\
x^2=\sinh{r} \cos{(\theta-\phi)}, ~ ~
x^3=\sinh{r} \sin{(\theta-\phi)},
\end{gathered} \tag{2.9}
\end{equation}
where $(r,\theta)$ are the polar coordinates of the base plane and $\phi$ is just the fibre coordinate. We note that
$$-x^0x^0-x^1x^1+x^2x^2+x^3x^3=-\cosh^2{r}+\sinh^2{r}=-1<0.$$
The inhomogeneous coordinates corresponding to (2.9), that play an important role in the later visualization of prism tilings in $\EUC$,
are given by
\begin{equation}
\begin{gathered}
x=\frac{x^1}{x^0}=\tan{\phi}, ~ ~
y=\frac{x^2}{x^0}=\tanh{r} \frac{\cos{(\theta-\phi)}}{\cos{\phi}}, \\
z=\frac{x^3}{x^0}=\tanh{r} \frac{\sin{(\theta-\phi)}}{\cos{\phi}}.
\end{gathered} \tag{2.10}
\end{equation}
\subsection{Geodesic curves}
The infinitesimal arc-length-square can be derived by the standard pull back method.
By $T^{-1}$-action of (2.6) on the differentials $(\mathrm{d}x^0;\mathrm{d}x^1;\mathrm{d}x^2;\mathrm{d}x^3)$, we obtain
that in this parametrization
the infinitesimal arc-length-square
at any point of $\SLR$ is the following:
\begin{equation}
   \begin{gathered}
      (\mathrm{d}s)^2=(\mathrm{d}r)^2+\cosh^2{r} \sinh^2{r}(\mathrm{d}\theta)^2+\big[(\mathrm{d}\phi)+\sinh^2{r}(\mathrm{d}\theta)\big]^2.
       \end{gathered} \tag{2.11}
     \end{equation}
Hence we get the symmetric metric tensor field $g_{ij}$ on $\SLR$ by components:
     \begin{equation}
       g_{ij}:=
       \begin{pmatrix}
         1&0&0 \\
         0&\sinh^2{r}(\sinh^2{r}+\cosh^2{r})& \sinh^2{r} \\
         0&\sinh^2{r}&1 \\
         \end{pmatrix}. \tag{2.12}
     \end{equation}

The geodesic curves of $\SLR$ are generally defined as having locally minimal arc length between any two of their (close enough) points.

By (2.9) the second order differential equation system of the $\SLR$ geodesic curve is the following:
\begin{equation}
\begin{gathered}
\ddot{r}=\sinh(2r)~\! \dot{\theta}~\! \dot{\phi}+\frac12 \big( \sinh(4r)-\sinh(2r) \big)\dot{\theta} ~\! \dot{\theta},\\
\ddot{\phi}=2\dot{r}\tanh{(r)}(2\sinh^2{(r)}~\! \dot{\theta}+ \dot{\phi}),\\ \ddot{\theta}=\frac{2\dot{r}}{\sinh{(2r)}}\big((3 \cosh{(2r)}-1)
\dot{\theta}+2\dot{\phi} \big). \tag{2.13}
\end{gathered}
\end{equation}

We can assume, by the homogeneity, that the starting point of a geodesic curve is the origin $(1,0,0,0)$.
Moreover, $r(0)=0,~ \phi(0)=0,~ \theta(0)=0,~ \dot{r}(0)=\cos(\alpha),~ \dot{\phi}(0)=\sin(\alpha)=-\dot{\theta}(0)$ are the initial values
in Table 1 for the solution of (2.12),
and so the unit velocity will be achieved.
\smallbreak
\centerline{\vbox{
\halign{\strut\vrule~\hfil $#$ \hfil~\vrule
&\quad \hfil $#$ \hfil~\vrule
&\quad \hfil $#$ \hfil\quad\vrule
&\quad \hfil $#$ \hfil\quad\vrule
&\quad \hfil $#$ \hfil\quad\vrule
\cr
\noalign{\vskip2pt}
\noalign{\hrule}
\multispan2{\strut\vrule\hfill{\bf  Table 1} \hfill\vrule}%
\cr
\noalign{\hrule}
\noalign{\vskip2pt}
{\rm Types} & {}  \cr
\noalign{\vskip2pt}
\noalign{\hrule}
\noalign{\vskip2pt}
\begin{gathered} 0 \le \alpha < \frac{\pi}{4} \\ (\bH^2-{\rm like~direction}) \end{gathered}
& \begin{gathered}  r(s,\alpha)={\mathrm{arsinh}} \Big( \frac{\cos{\alpha}}{\sqrt{\cos{2\alpha}}}\sinh(s\sqrt{\cos{2\alpha}}) \Big) \\
\theta(s,\alpha)=-{\mathrm{arctan}} \Big( \frac{\sin{\alpha}}{\sqrt{\cos{2\alpha}}}\tanh(s\sqrt{\cos{2\alpha}}) \Big) \\
\phi(s,\alpha)=2\sin{\alpha} s + \theta(s,\alpha) \end{gathered} \cr
\noalign{\vskip2pt}
\noalign{\hrule}
\noalign{\vskip2pt}
\begin{gathered} \alpha=\frac{\pi}{4} \\ ({\rm light~direction}) \end{gathered} &
\begin{gathered}  r(s,\alpha)={\mathrm{arsinh}} \Big( \frac{\sqrt{2}}{2} s \Big) \\
\theta(s,\alpha)=-{\mathrm{arctan}} \Big( \frac{\sqrt{2}}{2} s \Big) \\
\phi(s,\alpha)=\sqrt{2} s +\theta(s,\alpha) \end{gathered} \cr
\noalign{\vskip1pt}
\noalign{\hrule}
\noalign{\vskip2pt}
\begin{gathered} \frac{\pi}{4}  < \alpha \le \frac{\pi}{2} \\ ({\rm fibre-like~direction}) \end{gathered} &
\begin{gathered}  r(s,\alpha)={\mathrm{arsinh}} \Big( \frac{\cos{\alpha}}{\sqrt{-\cos{2\alpha}}}\sin(s\sqrt{-\cos{2\alpha}}) \Big) \\
\theta(s,\alpha)=-{\mathrm{arctan}} \Big( \frac{\sin{\alpha}}{\sqrt{-\cos{2\alpha}}}\tan(s\sqrt{-\cos{2\alpha}}) \Big) \\
\phi(s,\alpha)=2\sin{\alpha} s + \theta(s,\alpha) \end{gathered}  \cr
\noalign{\vskip2pt}
\noalign{\hrule}
\noalign{\hrule}}}}
\smallbreak
The equation of the geodesic curve in the hyperboloid model can be given by the usual geographical sphere coordinates
$(\lambda, \alpha)$, as longitude and altitude, respectively, from the general starting position of (2.10), (2.11),
$(-\pi < \lambda \le \pi, ~ -\frac{\pi}{2}\le \alpha \le \frac{\pi}{2})$,
and the arc-length parameter $0 \le s \in \bR$. The Euclidean coordinates $X(s,\lambda,\alpha)$,
$Y(s,\lambda,\alpha)$, $Z(s,\lambda,\alpha)$ of the geodesic curves can be determined by substituting the results of Table 1 into the
equations (2.11) as follows (see e.g. \cite{MSz})
\begin{equation}
\begin{gathered}
X(s,\lambda,\alpha)=\tan{(\phi(s,\alpha))}, \\
Y(s,\lambda,\alpha)=\frac{\tanh{(r(s,\alpha))}}{\cos{(\phi(s,\alpha)}} \cos \big[ \theta(s,\alpha)-\phi(s,\alpha)+\lambda \big],\\
Z(s,\lambda,\alpha)=\frac{\tanh{(r(s,\alpha))}}{\cos{(\phi(s,\alpha)}} \sin \big[ \theta(s,\alpha)-\phi(s,\alpha)+\lambda \big].
\end{gathered} \tag{2.14}
\end{equation}
\begin{Definition}
The {\rm geodesic distance} $d(P_1,P_2)$ between the points $P_1$ and $P_2$ is defined by the arc length of the geodesic curve
from $P_1$ to $P_2$.
\end{Definition}
\section{On regular prism tilings}
In \cite{Sz13-1} we defined and described the regular prisms and prism tilings with a space group class $\Gamma=\mathbf{pq2_1}$ of $\SLR$.
These will be summarized in this section.
\begin{Definition}[\cite{Sz13-1}]
Let $\cP^i$ be an infinite solid that is bounded by certain surfaces
determined (as in \cite{Sz13-1}) by ``side fibre lines" passing through the
vertices of a regular $p$-gon $\cP^b$ lying in the base plane.
The images of solids $\cP^i$ by $\SLR$ isometries are called {\rm infinite regular $p$-sided prisms}.
Here regular means that the side surfaces are congruent to each other under rotations about a fiber
line (e.g. through the origin).
\end{Definition}
The common part of $\cP^i$ with the base plane is the {\it base figure} of $\cP^i$ that is denoted by $\cP(p,q)$ and its vertices coincide
with the vertices of $\cP^b$, {\bf but $\cP(p,q)$ is not assumed to be a polygon}.
\begin{Definition}[\cite{Sz13-1}]
A {\rm bounded regular $p$-sided prism} is an isometric image of a solid
which is bounded by the side surfaces of a regular $p$-sided infinite prism $\cP^i$ its base figure
$\cP(p,q)$ and the translated copy $\cP^t(p,q)$ of $\cP(p,q)$ by a fibre translation, given by (2.2).
The faces $\cP(p,q)$ and $\cP^t(p,q)$ are called {\rm cover faces}.
\end{Definition}
We consider regular prism tilings $\cT_p(q)$ by prisms $\cP_p(q)$ where $q$ pieces regularly meet
at each side edge by $q$-rotation.

The following theorem has been proved in \cite{Sz13-1}:
\begin{theorem}[\cite{Sz13-1}]
There exist regular bounded not face-to-face prism tilings $\cT_p(q)$ in $\SLR$ for all integers $p$ and $q$ such that $3 \le p $ and
$\frac{2p}{p-2} < q$.
\end{theorem}
We assume that the prism $\cP_p(q)$ is a {\it topological polyhedron} having at each vertex
one $p$-gonal cover face (it is not a polygon at all) and two {\it skew quadrangles} which lie on certain side surfaces in the model.
Let $\cP_p(q)$ be one of the tiles of $\cT_p(q)$, $\cP^b$ is centered in the origin with vertices $A_1A_2A_3 \dots A_p$ in the base plane (Fig.~1 and 2).
It is clear that the side curves $c_{A_iA_{i+1}}$ $(i=1\dots p, ~ A_{p+1} \equiv A_1)$
of the base figure are derived from each other by $\frac{2\pi}{p}$ rotation about the vertical $x$ axis, so there are congruent in $\SLR$ sense.
The corresponding vertices $B_1B_2B_3 \dots B_p$ are generated by a fibre translation $\tau$ given by (2.3)
with real parameter $\Phi > 0$.
\begin{figure}[ht]
\centering
\includegraphics[width=12cm]{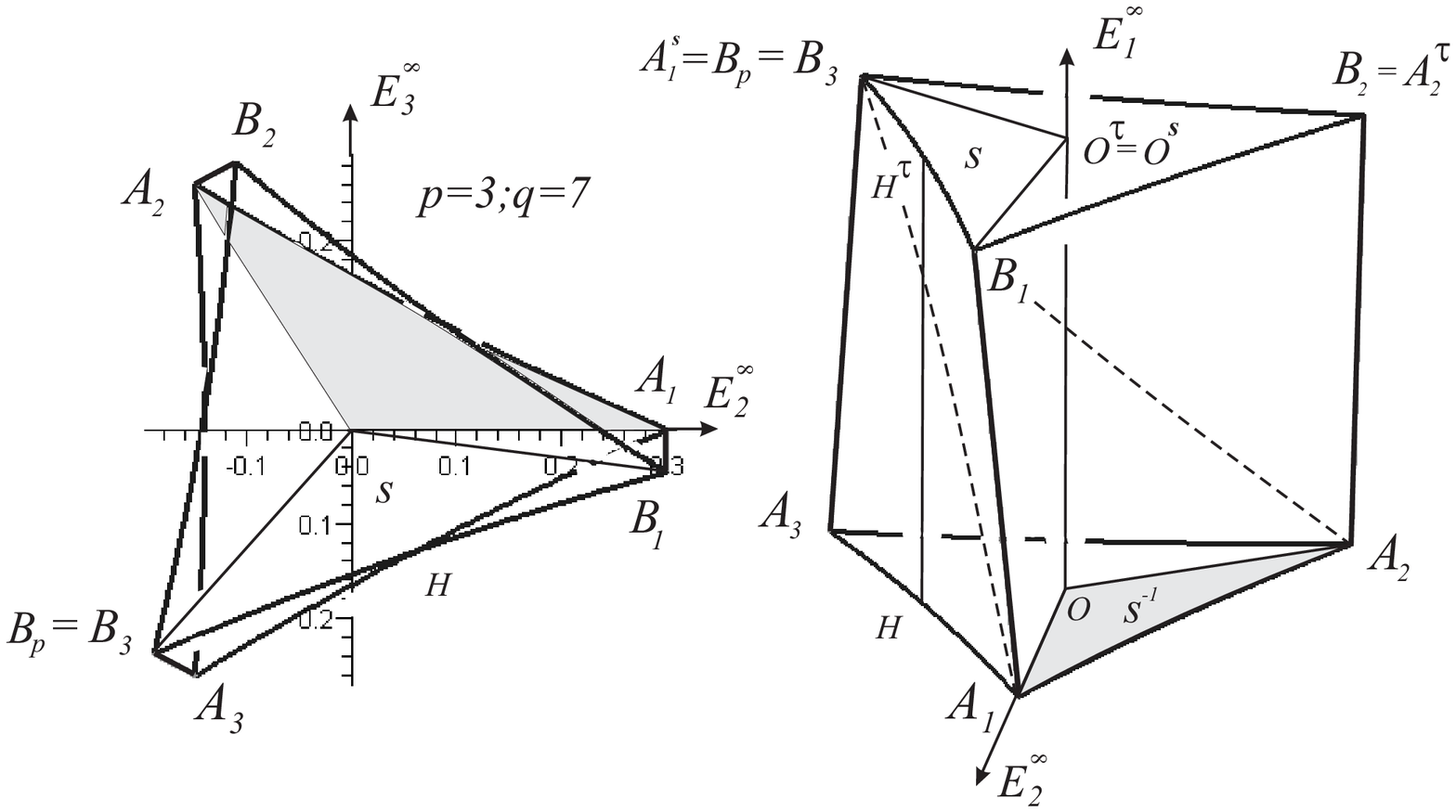}
\caption{The regular prism $\cP_p(q)$ and the fundamental domain of the space group ${\mathbf{pq2_1}}$}
\label{}
\end{figure}
The fibre lines through the vertices
$A_iB_i$ are denoted by $f_i, \ (i=1, \dots, p)$ and the fibre line through the ``midpoint" $H$ of the curve $c_{A_1A_{p}}$ is denoted by
$f_0$. This $f_0$ will be a half-screw axis as follows below.

The tiling $\cT_p(q)$ is generated by a
discrete isometry group $\Gamma_p(q)=\mathbf{pq2_1}$ $\subset Isom(\SLR)$
given by its fundamental domain $A_1A_2O A_1^{\bs} A_2^{\bs} O^{\bs}$ a {\it topological polyhedron} and the group presentation
(see Fig.~1 and 4 for $p=3$ and \cite{Sz13-1} for details):
\begin{equation}
\begin{gathered}
\mathbf{pq2_1}=\{ \ba,\bb,\bs: \ba^p=\bb^q=\ba \bs \ba^{-1} \bs^{-1}= \bb \ba \bb \bs^{-1}=\mathbf{1} \}= \\
= \{ \ba,\bb: \ba^p=\bb^q=\ba \bb \ba \bb \ba^{-1} \bb^{-1} \ba^{-1} \bb^{-1}=\mathbf{1} \}. \tag{3.1}
\end{gathered}
\end{equation}
Here $\ba$ is a ${p}$-rotation about the fibre line through the origin ($x$ axis),  $\bb$ is a ${q}$-rotation about the fibre line trough
$A_1$ and $\bs=\bb \ba \bb$ is a screw motion~ $\bs:~ OA_1A_2 \rightarrow O^{\bs} B_p B_1$. All these can be obtained by formulas (2.5) and (2.6).
Then we get that $\ba\bb\ba\bb=\bb\ba\bb\ba=:\tau$ is a fibre translation.
Then $\ba \bb$ is a $\mathbf{2_1}$ half-screw motion about
$f_0=HH^{\tau}$ (see Fig.~1) that also determines the fibre translation $\tau$ above. This group in (3.1) surprisingly occurred in \S~6 of paper \cite{MSzV} at double links
$K_{p,q}$.
The coordinates of the vertices $A_1A_2A_3 \dots A_p$ of the base figure and the corresponding vertices $B_1B_2B_3 \dots B_p$
of the cover face can be computed for all given parameters $p,q$ by
\begin{equation}
\tanh(OA_1)=b:=\sqrt{\frac{1-\tan{\frac{\pi}{p}} \tan{\frac{\pi}{q}}} {1+\tan{\frac{\pi}{q}} \tan{\frac{\pi}{q}}}}. \tag{3.2}
\end{equation}
The volume formula of a {\it sector-like} 3-dimensional domain $Vol(D(\Psi))$ can be computed routinely by the metric tensor $g_{ij}$ (see \cite{MSz14})
in the hyperboloid coordinates.
This defined by the base figure $D$ lying in the base plane and by
fibre translation $\tau$ given by (2.3) with the height parameter $\Psi$.
\begin{theorem}[\cite{MSz14}]
Suppose we are given a sector-like region $D$, so a continuous function $r = r(\theta)$
where the radius $r$ depends upon the polar angle $\theta$. The volume of the domain $D(\Psi)$ is derived by the following integral:
\begin{equation}
\begin{gathered}
Vol(D(\Psi))=\int_{D}  \frac{1}{2}\sinh(2r(\theta)) {\mathrm{d}}r~ {\mathrm{d}} \theta  ~{\mathrm{d}} \psi = \\ = \int_0^{\Psi} \int_{\theta_1}^{\theta_2} \int_{0}^{r(\theta)}
\frac{1}{2}\sinh(2r(\theta))~ {\mathrm{d}}r~ {\mathrm{d}} \theta \ {\mathrm{d}} \psi =\Psi \int_{\theta_1}^{\theta_2}
\frac{1}{4}(\cosh(2r(\theta))-1)~ {\mathrm{d}} \theta.
\tag{3.3}
\end{gathered}
\end{equation}
\end{theorem}
Let $\cP_p(q)$ be an arbitrary bounded regular prism. We get the following
\begin{theorem}[\cite{MSz14}]
The volume of the bounded regular prism $\cP_p(q)$ \Big($3 \le p \in \mathbb{N} $, $\frac{2p}{p-2} < q \in \mathbb{N}$\Big) is given by the following simple formula:
\begin{equation}
{\mathrm{Vol}}(\cP_p(q))={\mathrm{Vol}}(D(p,q,\Psi)) \cdot p, \tag{3.4}
\end{equation}
where ${\mathrm{Vol}}(D(p,q,\Psi))$ is the volume of the sector-like 3-dimensional domain given by the sector region $OA_1A_2 \subset \cP(p.q)$ (see Fig.~1. and Fig.~3.)
and by $\Psi$, the $\SLR$ height of the prism, depending on $p,q$.
\end{theorem}
\begin{corollary}
It is clear, that if $\Psi=1$ then we obtain the area of the base figur $\cP(p,q)$: ${\mathrm{Vol}}(\cP_p(q))=Vol(D(p,q,\Psi=1)) \cdot p$ $={\mathrm{Area}}(\cP(p,q))$.
\end{corollary}
\begin{figure}[ht]
\centering
\includegraphics[width=5cm]{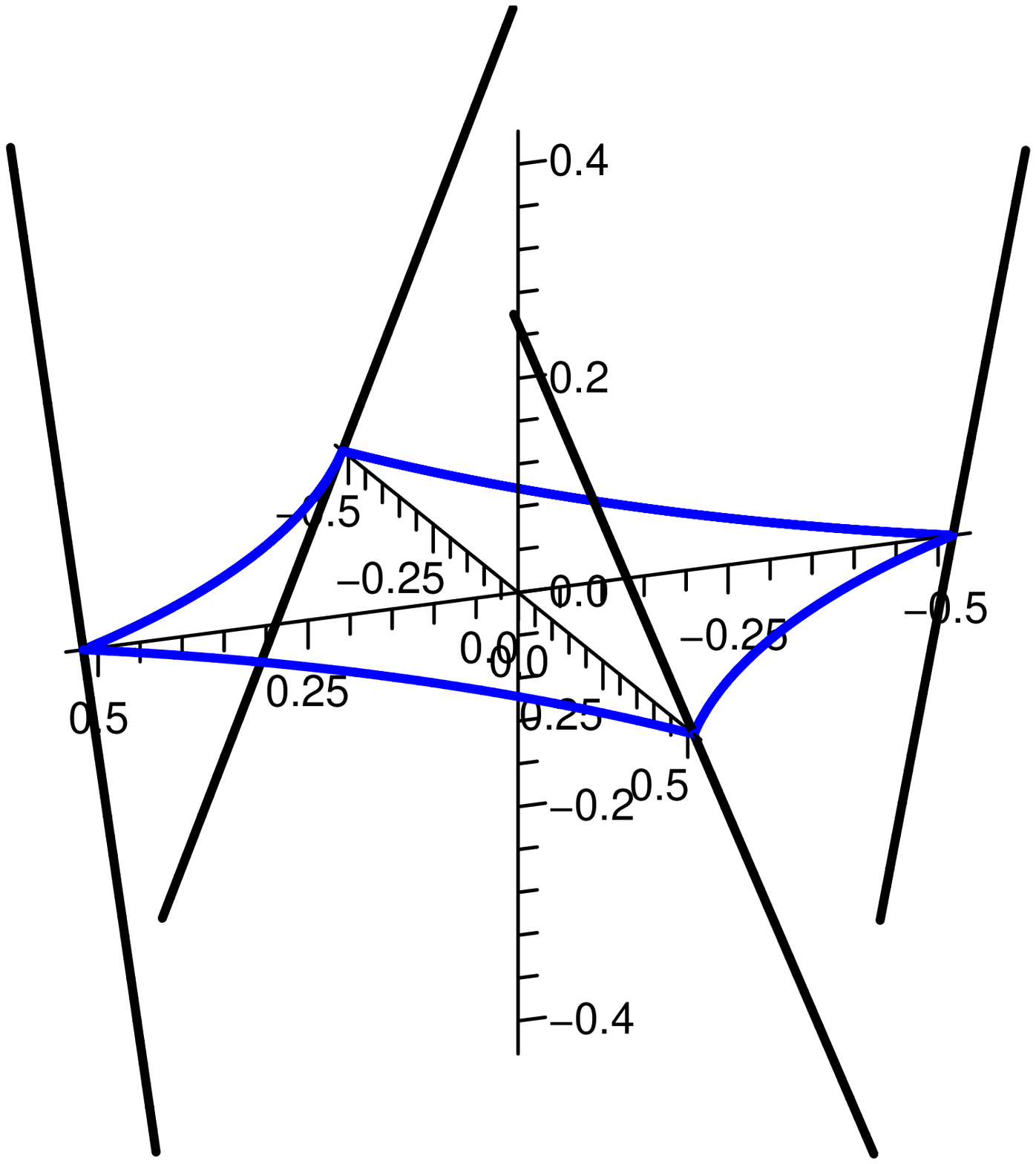} \includegraphics[width=6cm]{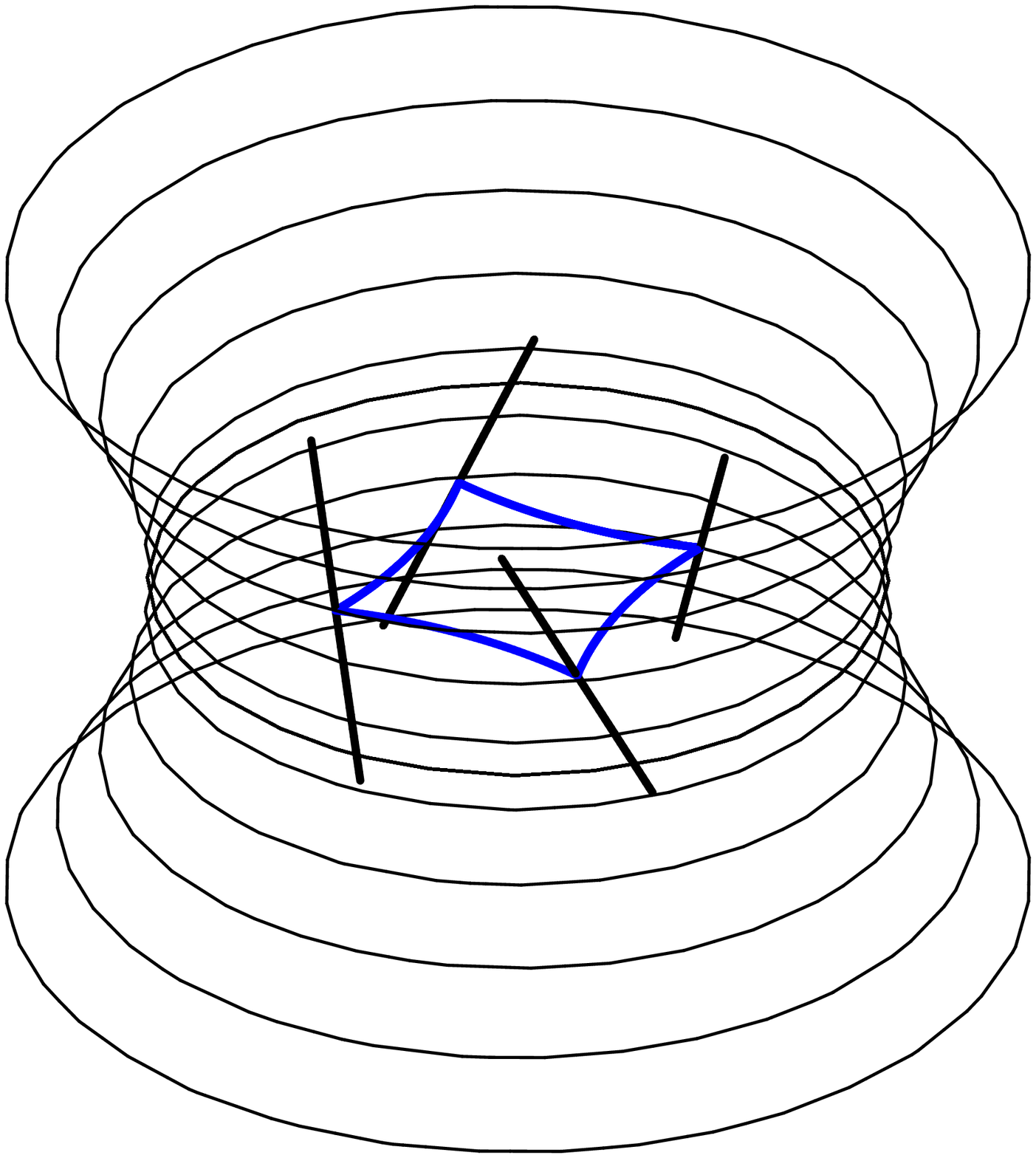}
\caption{Regular infinite 4-gonal prism $\cP_4^i(6)$ of the infinite regular prism tiling $\cT_4^i(6)$}
\label{}
\end{figure}
In this subsection we study the regular infinite prism tilings $\cT^i_p(q)$.
Let $\cT_p(q)$ be a regular prism tiling and let
$\cP_p(q)$ be one of its tiles given by its base figure $\cP$ centred at the origin $E_0$ with vertices $A_1A_2A_3 \dots A_p$
in the base plane of the model and the corresponding vertices $A_1A_2$ $A_3 \dots A_p$ and $B_1B_2B_3 \dots B_p$
are generated by fibre translation $\tau$ given by (2.3)
with parameter $\Psi=2\cdot(\frac{\pi}{2}-\frac{\pi}{p}-\frac{\pi}{q})$. The images of the topological polyhedron $\cP_p(q)$
by the translations $\langle \tau \rangle $ form an infinite prism $\cP^i_p(q)$ (see Definitions 3.~1-2).
\begin{figure}[ht]
\centering
\includegraphics[width=12cm]{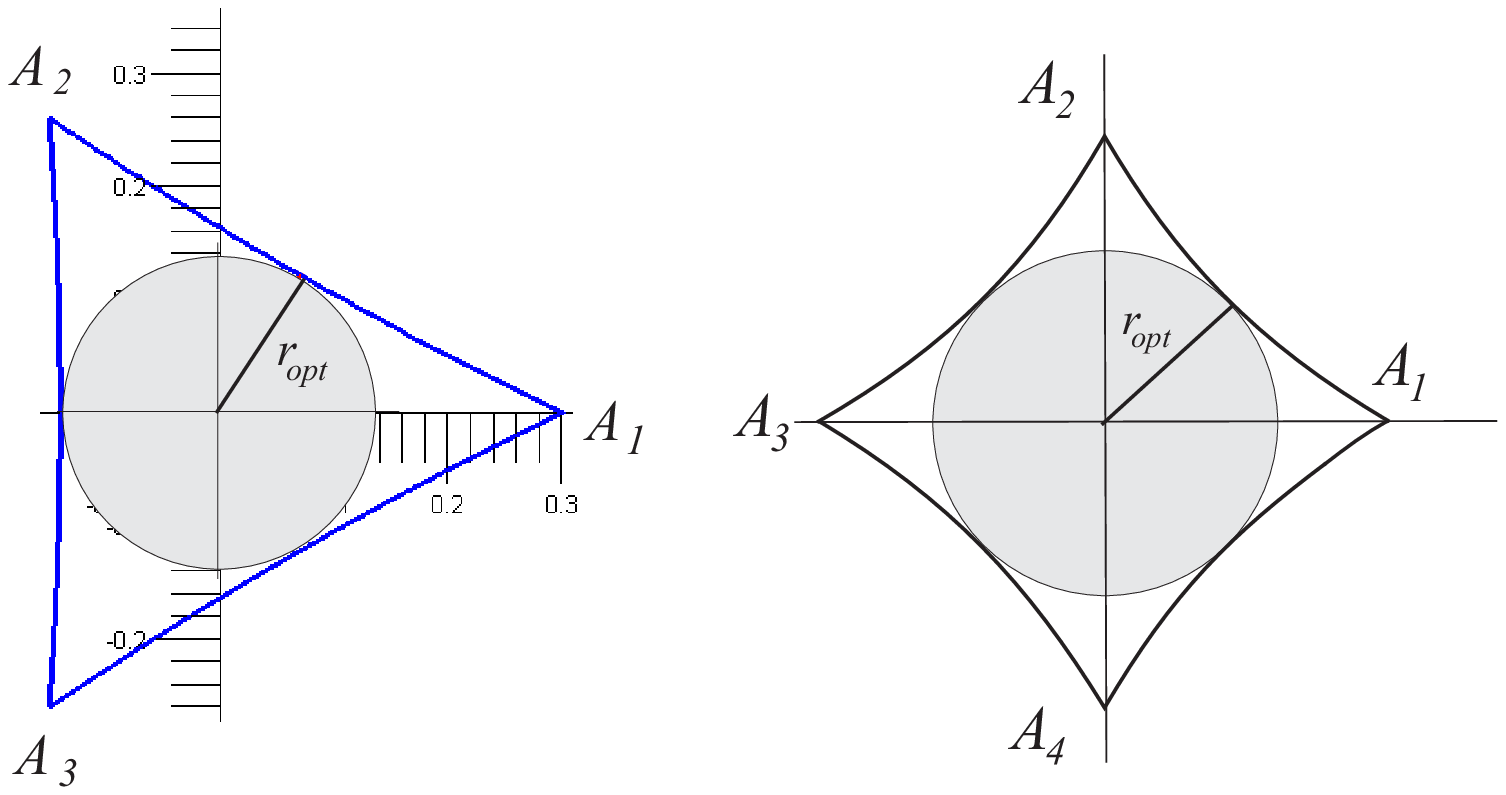}
\caption{The maximum radius $r_{opt}$ for parameters $p=3$, $q=7$ (left) and $p=4$, $q=6$ (right)}
\label{}
\end{figure}
By the construction of the bounded prism tilings it follows that the rotation through $\omega=\frac{2\pi}{q}$ about the fibre lines
$f_i$ maps the corresponding side face onto the neighbouring one.
Therefore, we have got the following (see \cite{Sz13-1}):
\begin{theorem}[\cite{Sz13-1}]
There exist regular infinite face-to-face prism tilings $\cT_p^i(q)$ for integer parameters $p\ge 3$ and $q > \frac{2p}{p-2}$.
\end{theorem}
For example, we have described $\cP_4^i(6)$ with its base polygon in Fig.~2, with $b=\frac{\sqrt{6}-\sqrt{2}}{2}$.
\section{Fibre-like cylinders, their packings and coverings}
\subsection{On $\SLR$ cylinders}
\begin{Definition}
Let $\cC^i(r)$ be an infinite solid that is bounded by certain surfaces
determined by ``side fibre lines" passing through the points of a circle $\cC^b(r)$ of radius $r\in \mathbf{R}^+$ lying in the base plane and centred at the origin.
The images $\cC^i_\bt(r)$ of solid $\cC^i(r)$ by $\SLR$ isometries $\bt$ are called {\rm infinite circular cylinders}.
\end{Definition}
The common part of $\cC^i_\bt(r)$ with the base plane is the {\it base figure} of $\cC^i_\bt(r)$ that is denoted by $\cC_\bt(r)$.
\begin{Definition}
A {\rm bounded fibre-like circular cylinder} is an isometric image of a solid
which is bounded by the side surface of a infinite circular cylinder $\cC^i(r)$,
its base figure
$\cC^b(r)$ and the translated copy $\cC^c(r)$ of $\cC^b(r)$ by a fibre translation, given by (2.2).
The faces $\cC^b(r)$ and $\cC^c(r)$ are called {\rm cover faces}.
The height (or altitude) of the cylinder is the distance between its cover faces.
\end{Definition}
The direct consequence of the (2.2) and (2.3) formulas and the Definitions 4.1-2 the following 
\begin{lemma}
If the radius of the circle $\cC^b(r)$ centred at the origin lying in the base plane is $r$ then the corresponding infinite fibre-like cylinder $\cC^i(r)$
is a one-sheeted hyperboloid in Euclidean sense with equation
\begin{equation}
\frac{y^2}{\tanh^2{r}}+\frac{z^2}{\tanh^2{r}}-x^2=1, \quad {\text{where}} \quad x=\frac{x^1}{x^0},~y=\frac{x^2}{x^0},~z=\frac{x^3}{x^0}. \tag{4.1}
\end{equation}
\end{lemma}
Let us denote the image of the cylinder $\cC^i(r)$ at the $\SLR$ translation $\bT$ (see (2.5))
by $\cC^i_\bT(r)$ where the translation is given by parameters $x^0=1,x^1,x^2,x^3$ and denote
its common part with the base plane by $\cC_\bT(r)$.

Using Lemma 1, formulas (2.5), (2.8) and the results of classical differential geometry, we obtain the following
\begin{lemma}
$\cC_\bT(r)$ is a circle of radius
\begin{equation}
R=\frac{\tanh^2{r}(1+(x^1)^2-(x^2)^2-(x^3)^2)}{1+(x^1)^2-\tanh^2{r}\cdot((x^2)^2+(x^3)^2)}, \tag{4.2}
\end{equation}
in Euclidean sense, centred at $K=(1;0;\frac{x^2-x^1x^3}{1+x^1x^1};\frac{x^3+x^1x^2}{1+x^1x^1})$. 
\end{lemma}
\begin{Remark}
\begin{enumerate}
\item The translations correspond a fibrum to a fibrum, so it is sufficient to consider the translations in the base plane, that means can be assumed $x^1=0$.
\item If a translation is given by $(1;0;x^2;x^3)$ where $\sqrt{(x^2)^2+(x^3)^2}=\tanh{2r}$ then the fibre-like cylinder $\cC^i(r)$ and its translated copies touch each other. 
In this case $R=\frac{\tanh{r}(1-\tanh^2{r})}{1-\tanh^2{2r}\tanh^2{r}}$ whose graph is shown in the Fig.~4. 
\end{enumerate}   
\end{Remark}
\begin{figure}[ht]
\centering
\includegraphics[width=6cm]{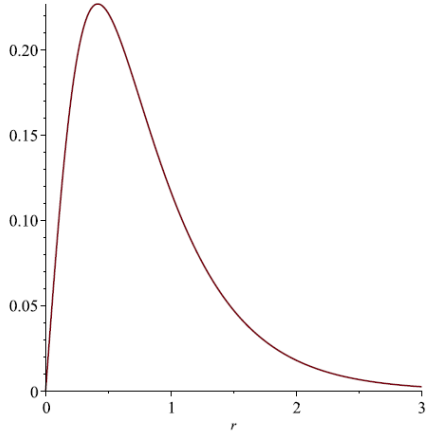}
\caption{The grapf of function $R=\frac{\tanh{r}(1-\tanh^2{r})}{1-\tanh^2{2r}\tanh^2{r}}$. It reaches its maximum in $r \approx 0.41572$.}
\label{}
\end{figure}
Using the infinitesimal arc-length-square (see (2.11)) and the corresponding metric tensor (2.12) we obtain the volume element in hyperboloid coordinates: 
\begin{equation}
\mathrm{d}V=\sqrt{\det(g_{ij})}~dr ~\mathrm{d} \theta ~\mathrm{d} \phi= \frac{1}{2}\sinh(2r) \mathrm{d}r ~\mathrm{d} \theta~ \mathrm{d} \phi.
\tag{4.3}
\end{equation}
Applying the (2.11), (2,12) and the previous formula, we directly obtain  the following lemmas:
\begin{lemma}
The perimeter $p(\cC(r))$ of a circle $\cC(r)$ with radius $r$ ($\SLR$ geodesic distance) lying in the base plane centred at the origin of the projective model of the $\SLR$
geometry can be calculated with the following formula
\begin{equation}
p(\cC(r))=\int_{-\pi}^\pi \sqrt{\sinh^2(r)(\sinh^2(r)+\cosh^2(R))}~\mathrm{d} \theta=2\pi\sinh(r)\sqrt{\cosh{(2r)}}. \tag{4.4}
\end{equation}
\end{lemma}
\begin{lemma}
The area ${\mathrm{Area}}(\cC(r))$ of a circle with radius $r$ lying in the base plane centred at the origin of the projective model of the $\SLR$
geometry can be calculated with the following formula
\begin{equation}
{\mathrm{Area}}(\cC(r))=\pi \cdot \int_0^r \sinh(2r) ~\mathrm{d}t =\pi \cdot \sinh^2(r).\tag{4.5}
\end{equation}
\end{lemma}
\begin{Corollary}
The surface area ${\mathrm{SA}}(\cC(r,\Psi))$ of a bounded fibre-like circular cylinder $\cC(r,\Psi)$ \Big($3 \le p \in \mathbf{N} $,
$\frac{2p}{p-2} < q \in \mathbf{N}$, $\Psi\in \mathbf{R}^+$ \Big)
is
\begin{equation}
{\mathrm{SA}}(\cC(r,\Psi))=2 \cdot \pi \sinh(r)(\sqrt{\cosh{(2r)}}\cdot \Psi+\sinh(r)). \notag
\end{equation}
\end{Corollary}

\begin{lemma}
The volume ${\mathrm{Vol}}((\cC(r,\Psi)))$ a bounded fibre-like circular cylinder with radius $r$ and height (or altitude) $\Psi$ in the projective model of the $\SLR$
geometry can be calculated with the following formula
\begin{equation}
{\mathrm{Vol}}((\cC(r,\Psi)))=\Psi \cdot \pi \cdot \sinh^2(r).\tag{4.6}
\end{equation}
\end{lemma}
\subsubsection{Optimal inscribed circle radius and packing density}
We consider an infinite regular prism tiling $\cT_p^i(q)$ and let $\cP_p^i(q)$ be one of its tiles with base figure $\cP(p,q)$
centered at the origin with vertices $A_1A_2 \dots A_p$ in the base plane of the model.
Let $\cC^i(r_{opt}(p,q))$ be a fibre-like infinite circular cylinder of radius $r_{opt}(p,q)$ with circular base figur $\cC^b(r_{opt}(p,q))$ that centred at the origin where
it is the inscribed circle of $\cP(p,q)$, that means, that $\cC^i(r_{opt}(p,q))$ cylinder touch the side surfaces of the prism $\cP_p^i(q)$.

It is obvious that if we consider the corresponding images of the cylinder $\cC^i(r_{opt}(p,q))$ to each element of the prism tiling $\cT_p^i(q)$, 
we get a cylinder packing arrangement in the $\SLR$ space. In these arrangements, neighbouring cylinders touch each other and one cylinder is touched by exactly $p$ others.

For the density of the optimal cylinder packing it is sufficient to relate the volume of the optimal (touching) infinite cylinder $\cC^i(r_{opt}(p,q))$
to that of the infinite prism $\cP_p^i(q)$. This ratio can be replaced by the ratio of the areas of base figurs ${\mathrm{Area}}(\cC^b(r_{opt}(p,q)))$ and ${\mathrm{Area}}(\cP(p,q))$
(see (3.3), (3.4) and (4.5)).
\begin{Definition}
The density of the optimal cylinder packing related to the prism tiling $\cT^{i}_p(q)$ (with integer parameters
$p \ge 3$ and $q > \frac{2p}{p-2}$) is the following:
\begin{equation}
\delta^{opt}_p(q):=\frac{{\mathrm{Area}}(\cC^b(r_{opt}(p,q)))}{{\mathrm{Area}}(\cP(p,q))} \notag
\end{equation}
where $r_{opt}(p,q))$ is the radius of the inscribed circle of $\cP(p,q)$.
\end{Definition}
In \cite{Sz14} we determined the equation of the side curve $c_{A_1A_2}$  as the foot points
(see Fig.~3) of the corresponding fibre lines ($3 \le p,~ \frac{2p}{p-2}<q$, where $p$ and $q$ are integer parameters):
\begin{lemma}[\cite{Sz14}]
The parametric equation of the side curve $c_{A_1A_2}$ of the base figur $\cP(p,q)$ is
\begin{equation}
\scriptsize
\begin{gathered}
c_p^q(t)=\Bigg(0,~\sqrt{\sin \left( \frac{2\pi}{p}+\frac{2\pi}{q} \right)}\Bigg( t \cos \left( \frac{2\pi}{p} \right)\sin^2\left(\frac{\pi}{p}+\frac{\pi}{q}\right)-
\frac{t}{2}\sin\left(\frac{2\pi}{p}\right)\sin\left(\frac{2\pi}{p}+\frac{2\pi}{q}\right)+ \\
\sin^2\left(\frac{\pi}{p}+\frac{\pi}{q}\right)(1-t)+ t^2 \cos\left(\frac{\pi}{p}+\frac{\pi}{q}\right)
\cos\left(\frac{\pi}{p}-\frac{\pi}{q}\right)\Bigg)\Big/ \\
\Bigg({\sqrt{\left( \sin \left( {\frac {2\pi }{p}} \right) +
\sin \left( {\frac {2\pi }{q}} \right)  \right)}} \Big(\sin^2\left(\frac{\pi}{p}+\frac{\pi}{q}\right)+t^2 \cos^2\left(\frac{\pi}{p}+\frac{\pi}{q}\right)\Big) \Bigg),\\
t\sqrt{\sin \left( \frac{2\pi}{p}+\frac{2\pi}{q} \right)}\Bigg( \sin \left( \frac{2\pi}{p} \right)\sin^2\left(\frac{\pi}{p}+\frac{\pi}{q}\right)+
\frac{1}{2}\cos\left(\frac{2\pi}{p}\right)\sin\left(\frac{2\pi}{p}+\frac{2\pi}{q}\right)(1-t)+ \\
\cos\left(\frac{\pi}{p}+\frac{\pi}{q}\right)\Big(t \sin\left(\frac{2\pi}{p}\right)\cos\left(\frac{\pi}{p}+\frac{\pi}{q}\right)+
\sin\left(\frac{\pi}{p}+\frac{\pi}{q}\right)(t-1)\Big)\Bigg)\Big/ \\
\Bigg({\sqrt{\left( \sin \left( {\frac {2\pi }{p}} \right) +
\sin \left( {\frac {2\pi }{q}} \right)  \right)}} \Big(\sin^2\left(\frac{\pi}{p}+\frac{\pi}{q}\right)+t^2
\cos^2\left(\frac{\pi}{p}+\frac{\pi}{q}\right)\Big) \Bigg), ~ ~ t \in [0,1].
\end{gathered} {\tag{\normalsize{4.7}}}
\end{equation}
\end{lemma}
The side curves $c_{A_iA_{i+1}}$ $(i=1\dots p, ~ A_{p+1} \equiv A_1)$
of the base figure are derived from each other by $\frac{2\pi}{p}$ rotation about the vertical $x$ axis, therefore they are congruent and
their curvatures are equal in $\SLR$ sense.
Moreover, the side curves are congruent also in the Euclidean sense, therefore their curvatures are equal in Euclidean sense, as well. 

We proved in \cite{Sz14} the following
\begin{lemma}[\cite{Sz14}]
The curvature $C_p(q)$ of the side curves $c_{A_iA_{i+1}}$ $(i=1\dots p, ~ A_{p+1} \equiv A_1)$ in the Euclidean sense is
\begin{equation}
C_p(q)=\sqrt {\frac{\cos \left( {\frac {\pi }{p}}+{\frac {\pi }{q}} \right)
\left( \sin \left( {\frac {2\pi }{p}} \right) +\sin \left( {
\frac {2\pi }{q}} \right)  \right)} { \sin \left( {\frac {\pi }{p}
}+{\frac {\pi }{q}}  \right)  \left( 1-\cos \left( {
\frac {2\pi }{p}} \right)  \right)}} \tag{4.8}
\end{equation}
therefore the side curves $c_{A_iA_{i+1}}$ $(i=1\dots p, ~ A_{p+1} \equiv A_1)$ are Euclidean circular arcs of radius $r_p^q=\frac{1}{C_p(q)}$.
\end{lemma}
From the above lemmas we obtain directly the radius $r_{opt}(p,q)$ of the inscribed circle of $\cP(p,q)$. 
\begin{lemma}
The radius $r_{opt}(p,q))$ of the inscribed circle of $\cP(p,q)$ is 
\begin{equation}
r_{opt}(p,q)=\mathrm{arctanh}\Bigg( \sqrt{\frac{\cos\frac{\pi}{q}-\sin\frac{\pi}{p}}{\cos\frac{\pi}{q}+\sin\frac{\pi}{p}}}\Bigg). \notag
\end{equation}
\end{lemma}
\subsubsection{On optimal circumscribed circle radius and covering density}
Similarly to the packing for the density of the optimal covering cylinder it is sufficient to relate the volume of the optimal circumscribing infinite cylinder $\cC^i(R_{opt}(p,q))$
to that of the infinite prism $\cP_p^i(q)$. This ratio can be replaced by the ratio of the areas of base figurs 
${\mathrm{Area}}(\cC^b(R_{opt}(p,q)))$ and ${\mathrm{Area}}(\cP(p,q))$ (see (3.3), (3.4) and (4.5)).
\begin{Definition}
The densitiy of the optimal cylinder covering related to the prism tiling $\cT^{i}_p(q)$ (with integer parameters
$p \ge 3$ and $q > \frac{2p}{p-2}$) is the following:
\begin{equation}
\Delta^{opt}_p(q):=\frac{{\mathrm{Area}}(\cC^b(R_{opt}(p,q)))}{{\mathrm{Area}}(\cP(p,q))}. \notag
\end{equation}
where $R_{opt}(p,q)$ is the radius of the circumscribed circle of $\cP(p,q)$.
\end{Definition}
It follows directly from formula (3.2) that 
\begin{equation}
R_{opt}(p,q)=E_0A_1=OA_1=\mathrm{arctanh}\Bigg({\sqrt{\frac{1-\tan{\frac{\pi}{p}} \tan{\frac{\pi}{q}}} {1+\tan{\frac{\pi}{q}} \tan{\frac{\pi}{q}}}}}\Bigg). \tag{4.9}
\end{equation}
\subsection{Optimal cylinder packing and covering densities}
We consider a infinite regular prism tiling $\cT_p^i(q)$. Let $\cP_p^i(q)$ be one of its tiles with base figure $\cP(p,q)$ and let 
$\cC^i(r_{opt}(p,q))$ be a fibre-like infinite circular cylinder of radius $r_{opt}$ with circular base 
figur $\cC^b(r_{opt}(p,q))$ that centred at the origin where it is the inscribed circle of $\cP(p,q)$. Moreover, let $\cC^i(R_{opt}(p,q))$ be the optimal circumscribing infinite 
cylinder of $\cP_p^i(q)$ with base figur $\cC^b(R_{opt}(p,q))$.

The radii $r_{opt}(p,q)$ and $R_{opt}(p,q)$ can be determined by applying the Lemma 4.8 and formula (4.9)
for all possible parameters. The areas ${\mathrm{Area}}(\cC^b(r_{opt}(p,q)))$ and ${\mathrm{Area}}(\cC^b(R_{opt}(p,q)))$ can be determined by Lemma 4.1. 
Moreover, the area of base figure $\cP(p,q)$ can be computed by Theorem 3.3.

The above locally densest geodesic cylinder packings and coverings can be determined for all regular prism tilings $\cT^i_p(q)$ ($p,q$ as above);
our results are summarized in Tables 2 and 3.
\medbreak
\centerline{\vbox{
\halign{\strut\vrule\quad \hfil $#$ \hfil\quad\vrule
&\quad \hfil $#$ \hfil\quad\vrule &\quad \hfil $#$ \hfil\quad\vrule &\quad \hfil $#$ \hfil\quad\vrule &\quad \hfil $#$ \hfil\quad\vrule
\cr
\noalign{\hrule}
\multispan5{\strut\vrule\hfill{Table 2, Packing} \hfill\vrule}%
\cr
\noalign{\hrule}
\noalign{\vskip2pt}
\noalign{\hrule}
(p, q)& r_{opt}(p,q) &  \mathrm{Area}(\cC^b(r_{opt}(p,q)) & \mathrm{Area}(\cP_p(q)) & \delta^{opt}_p(q) \cr
\noalign{\hrule}
\noalign{\hrule}
(3,7) &  0.14156 & 0.06338 & 0.11220 & 0.56489 \cr
\noalign{\hrule}
(3,8) &  0.18176 & 0.10494 & 0.19635 & 0.53443 \cr
\noalign{\hrule}
(3,10) & 0.21980 & 0.15423 & 0.31416 & 0.49093 \cr
\noalign{\hrule}
(3,1000) & 0.27465 & 0.24299 & 0.78069 & 0.31126 \cr
\noalign{\hrule}
\vdots & \vdots & \vdots & \vdots & \vdots \cr
\noalign{\hrule}
(4,5) &  0.26532 & 0.22639 & 0.31416 & 0.72061 \cr
\noalign{\hrule}
(4,6) &  0.32924 & 0.35303 & 0.52360 & 0.67424 \cr
\noalign{\hrule}
(4,10) & 0.40423 & 0.54192 & 0.94248 & 0.57500 \cr
\noalign{\hrule}
(4,1000) & 0.44068 & 0.65063 & 1.56451 & 0.41587 \cr
\noalign{\hrule}
\vdots & \vdots & \vdots & \vdots & \vdots \cr
\noalign{\hrule}
(7,3) &  0.27264 & 0.23936 & 0.26180 & 0.91430 \cr
\noalign{\hrule}
(7,4) &  0.53520 & 0.98915 & 1.17810 & 0.83962 \cr
\noalign{\hrule}
(7,5) &  0.61750 & 1.35810 & 1.72788 & 0.78600 \cr
\noalign{\hrule}
(7,10) & 0.71065 & 1.87233 & 2.82743 & 0.66220 \cr
\noalign{\hrule}
(7,1000) & 0.73867 & 2.04950 & 3.91600 & 0.52337 \cr
\noalign{\hrule}
\vdots & \vdots & \vdots & \vdots & \vdots \cr
\noalign{\hrule}
(9,3) &  0.46377 & 0.725552 & 0.78540 & 0.92380 \cr
\noalign{\hrule}
\vdots & \vdots & \vdots & \vdots & \vdots \cr
\noalign{\hrule}
(10,3) &  0.53064 & 0.97081 & 1.04720 & 0.92705 \cr
\noalign{\hrule}
\vdots & \vdots & \vdots & \vdots & \vdots \cr
\noalign{\hrule}
(20,3) &  0.91485 & 3.44983 & 3.66519 & 0.94124 \cr
\noalign{\hrule}
(40,3) &  1.26948 & ?????? & 8.90118 & 0.94813 \cr
\noalign{\hrule}
\vdots & \vdots & \vdots & \vdots & \vdots \cr
\noalign{\hrule}
(100,3) &  1.72981 & 23.43332 & 24.60914 & 0.95222 \cr
\noalign{\hrule}
\vdots & \vdots & \vdots & \vdots & \vdots \cr
\noalign{\hrule}
(1000,3) &  2.88151 & 248.42962 & 260.22859 & 0.95466 \cr
\noalign{\hrule}
(5000,3) &  3.68623 & 1248.429286 & 1307.42614 & 0.95488 \cr
\noalign{\hrule}
(p\rightarrow \infty,3) &   &  &  & \rightarrow 3/\pi \cr
\noalign{\hrule}
}}}
\medbreak
\medbreak
\centerline{\vbox{
\halign{\strut\vrule\quad \hfil $#$ \hfil\quad\vrule
&\quad \hfil $#$ \hfil\quad\vrule &\quad \hfil $#$ \hfil\quad\vrule &\quad \hfil $#$ \hfil\quad\vrule &\quad \hfil $#$ \hfil\quad\vrule
\cr
\noalign{\hrule}
\multispan5{\strut\vrule\hfill{Table 3, Covering} \hfill\vrule}%
\cr
\noalign{\hrule}
\noalign{\vskip2pt}
\noalign{\hrule}
(p, q)& R_{opt}(p,q) &  \mathrm{Area}(\cC^b(R_{opt}(p,q)) & \mathrm{Area}(\cP_p(q)) & \Delta^{opt}_p(q) \cr
\noalign{\hrule}
\noalign{\hrule}
(3,7) &  0.31034 & 0.31240 & 0.11220 & 2.78432 \cr
\noalign{\hrule}
(3,8) &  0.43035 & 0.61865 & 0.19635 & 3.15078 \cr
\noalign{\hrule}
(3,10) & 0.58867 & 1.22035 & 0.31416 & 3.88451 \cr
\noalign{\hrule}
(3,1000) & 2.95343 & 287.10339 & 0.78069 & 367.75794 \cr
\noalign{\hrule}
\vdots & \vdots & \vdots & \vdots & \vdots \cr
\noalign{\hrule}
(4,5) &  0.42124 & 0.59122 & 0.31416 & 1.88191 \cr
\noalign{\hrule}
(4,6) &  0.57311 & 1.14990 & 0.52360 & 2.19615 \cr
\noalign{\hrule}
(4,10) & 0.89491 & 3.26362 & 0.94248 & 3.46281 \cr
\noalign{\hrule}
(4,1000) & 3.22808 & 498.42756 & 1.56451 & 318.58317 \cr
\noalign{\hrule}
\vdots & \vdots & \vdots & \vdots & \vdots \cr
\noalign{\hrule}
(7,3) &  0.31034 & 0.31240 & 0.26180 & 1.19328 \cr
\noalign{\hrule}
(7,4) &  0.68003 & 1.69100 & 1.17810 & 1.43536 \cr
\noalign{\hrule}
(7,5) &  0.85559 & 2.91868 & 1.72788 & 1.68917 \cr
\noalign{\hrule}
(7,10) & 1.27092 & 8.46797 & 2.82743 & 2.99493 \cr
\noalign{\hrule}
(7,1000) & 3.59343 & 1036.68649 & 3.91600 & 264.73129 \cr
\noalign{\hrule}
\vdots & \vdots & \vdots & \vdots & \vdots \cr
\noalign{\hrule}
(9,3) &  0.51794 & 0.92089 & 0.78540 & 1.17251 \cr
\noalign{\hrule}
\vdots & \vdots & \vdots & \vdots & \vdots \cr
\noalign{\hrule}
(10,3) &  0.58867 & 1.22035 & 1.04720 & 1.16535 \cr
\noalign{\hrule}
\vdots & \vdots & \vdots & \vdots & \vdots \cr
\noalign{\hrule}
(20,3) &  0.98360 & 4.15514 & 3.66519 & 1.13368 \cr
\noalign{\hrule}
(40,3) &  1.34063 & 9.95246 & 8.90118 & 1.11811 \cr
\noalign{\hrule}
\vdots & \vdots & \vdots & \vdots & \vdots \cr
\noalign{\hrule}
(100,3) &  1.80161 & 27.28722 & 24.60914 & 1.10883 \cr
\noalign{\hrule}
\vdots & \vdots & \vdots & \vdots & \vdots \cr
\noalign{\hrule}
(1000,3) &  2.95343 & 287.10339 & 260.22859 & 1.10327 \cr
\noalign{\hrule}
(5000,3) &  3.75815 & 1441.80469 & 1307.42614 & 1.10278 \cr
\noalign{\hrule}
(p\rightarrow \infty,3) &   &  &  & \rightarrow \sqrt{12}/\pi \cr
\noalign{\hrule}
}}}
\medbreak
\begin{figure}[ht]
\centering
\includegraphics[width=14cm]{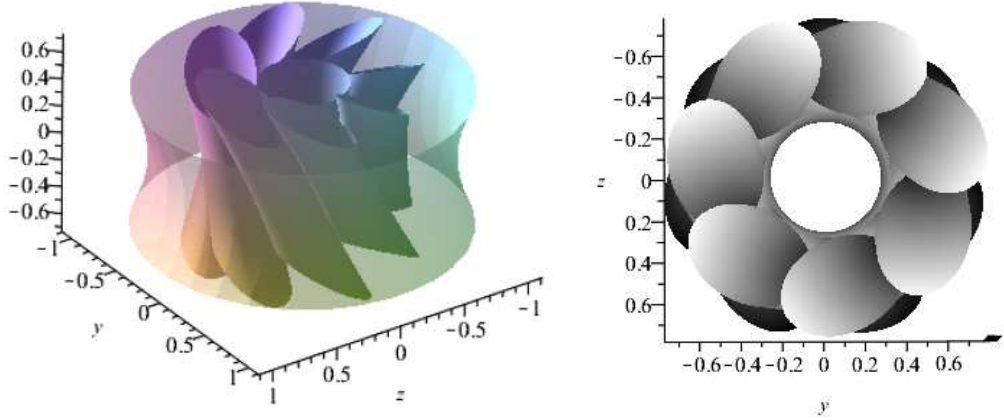}
\caption{The optimal cylinder packing configuration related to the prism tiling $\cT^{i}_7(3)$ in $\SLR$ space.}
\label{}
\end{figure}
\section{Connection with circle pacings and coverings belonging to regular hyperbolic mosaics}
It is well known, that there are infinitely many regular tessellations $\cM_p^q$ of the hyperbolic plane $\mathbf{H}^2$ with Schl\"afli symbol 
$\{p,q\}$ where the sum $\frac{1}{p} + \frac{1}{q} < \frac{1}{2},~(p,q \in \mathbf{N} ~ p,q \ge 3)$. Let $M_p^q$ be one of its tiles. Let 
$C(r^h(p,q))$ be its inscribed circle of radius $r^h(p,q)$ and $C(R^h(p,q))$ its circumscribed circle of radius $R^h(p,q)$. It is well known that
\begin{equation}
\begin{gathered}
r^h(p,q)={\mathrm{arccosh}} \Bigg( \frac{\cos\frac{\pi}{q}}{\sin\frac{\pi}{p}}\Bigg), \ \ R^h(p,q)={\mathrm{arccosh}}\Big(\cot\frac{\pi}{p}\cdot \cot\frac{\pi}{q}\Big) ~ 
{\mathrm{and}} \\ \tag{4.10}
\mathrm{Area}{(M_p^q)}=2\cdot p\cdot \Big(\frac{\pi}{2}-\frac{\pi}{p}-\frac{\pi}{q}\Big).
\end{gathered} 
\end{equation}
Examining the values obtained in (4.10) and the radii determined in Lemma 4.8, formula (4.9), furthermore, $\mathrm{Area}{(M_p^q)}$ and the area obtained from Theorem 3.2, we obtain the following
\begin{theorem}
\begin{equation}
\begin{gathered}
r^h(p,q)=2\cdot r_{opt}(p,q); \quad R^h(p,q)=2\cdot R_{opt}(p,q) \Rightarrow \\ 
{\mathrm{Area}}(C(r^h(p,q)))=4\cdot \pi \cdot \sinh^2\Big(\frac{r^h(p,q)}{2}\Big)= \\ 
=4 \cdot{\mathrm{Area}}(\cC(r_{opt}(p,q)))=4 \cdot \pi \cdot \sinh^2(r_{opt}(p,q)),\\
{\mathrm{Area}}(C(R^h(p,q)))=4\cdot \pi \cdot \sinh^2\Big(\frac{R^h(p,q)}{2}\Big)= \\ 
=4 \cdot{\mathrm{Area}}(\cC(R_{opt}(p,q)))=4 \cdot \pi \cdot \sinh^2(R_{opt}(p,q)), \\
\mathrm{Area}{(M_p^q)}=4\cdot \mathrm{Area}{(\cP)}.
\end{gathered} \notag
\end{equation}
\end{theorem}
\begin{Corollary}
\begin{enumerate}
\item The optimal cylinder packing density $\delta^{opt}_p(q)$ ~(see Definition 4.3)
in $\SLR$ space is equal of the densest hyperbolic circle packing density 
related to the regular mosaics with Schl\"afli symbol $\{p,q\}$.
\item The optimal cylinder covering density $\Delta^{opt}_p(q)$ ~ (see Definition 4.4) 
in $\SLR$ space is equal of the thinnest hyperbolic circle covering density 
related to the regular mosaics with Schl\"afli symbol $\{p,q\}$.
\end{enumerate}
\end{Corollary}
In a natural way similarly to the above $\SLR$ circular cylinder packings and 
coverings can be introduced circular cylinder packings in 
the hyperbolic space $\mathbf{H}^3$ where the axes of 
the circular cylinders are parallel, and in the $\HXR$ space, 
where the axes coincide with the fibrum lines of the space.

{\it Let us introduce the name fibrum type circular cylinders 
as a summary name for the previous circular cylinders.}

The following theorem is obtained directly from the previous results
\begin{theorem}
\begin{enumerate}
\item The optimal fibrum type cylinder packing density $\delta_c$ in $\HYP$, $\HXR$ and $\SLR$ spaces are equal of the densest hyperbolic circle packing density 
related to the regular mosaics with Schl\"afli symbol $\{p,q\}$, 
$\delta_c=\frac{3}{\pi}$.
\item The optimal fibrum type cylinder covering density $\Delta_c$ in $\HYP$, $\HXR$ and $\SLR$ spaces are equal of the densest hyperbolic circle covering density 
related to the regular mosaics with Schl\"afli symbol 
$\{p,q\}$, $\Delta_c=\frac{\sqrt{12}}{\pi}$.
\end{enumerate}
\end{theorem}

In this paper we mentioned only some natural problems related to 
Thurston spaces, but we hope that from these
the reader can appreciate that our projective method is suitable to study and solve 
similar problems that represent a huge class of open mathematical problems 
(see e.g. \cite{CsSz23}, \cite{MSz12}, \cite{M-Sz}, \cite{stachel}, \cite{Sz07}, \cite{Sz10-2}, \cite{Sz10-3}, \cite{Sz12-1}, 
\cite{Sz202}, \cite{Sz22}, \cite{Sz22-2}, \cite{Sz22-3}).
Detailed studies are the objective of ongoing research.


\begin{thebibliography}{MPSz98}
%
\bibitem{BF64} B\"or\"oczky,~K.~--~Florian,~A.
\"Uber die dichteste Kugelpackung im hyperbolischen Raum, \textit{Acta Math. Hung.}.
\bf{15} (1964), \rm, 237--245.
%
\bibitem{CsSz16}
Csima,~G.~--~Szirmai,~J.:
Interior angle sum of translation and geodesic triangles in $\SLR$ space, 
{\it Filomat}, {\bf 32/14} (2018), 5023--5036. 
%
\bibitem{CsSz23}
Csima,~G.~--~Szirmai,~J.:
Translation-like isoptic surfaces and angle sums of translation triangles in $\NIL$ geometry, 
{\it Submitted manuscript}, (2023). 
%
\bibitem{M97}
{Moln{\'a}r,~E.}
The projective interpretation of the eight 3-di\-men\-sional homogeneous geometries.
\emph{Beitr. Algebra Geom.,}
{\bf38} (1997), No.~2, 261--288.
  %
\bibitem{MSz}
{Moln{\'a}r,~E.~--~Szirmai,~J.}
Symmetries in the 8 homogeneous 3-geometries.
\emph{Symmetry Cult. Sci.,}
{\bf 21/1-3} (2010), 87-117.
%
\bibitem{MSz12}
{Moln{\'a}r,~E.~--~Szirmai,~J.}
Classification of $\SOL$ lattices.
\emph{Geom. Dedicata,}
{\bf 161/1} (2012), 251-275, DOI: 10.1007/s10711-012-9705-5.
%
\bibitem{MSz14}
{Moln{\'a}r,~E.~--~Szirmai,~J.}
Volumes and geodesic ball packings to the regular prism tilings in $\SLR$ space.
\emph{Publ. Math. Debrecen,}
{\bf 84/1-2} (2014), 189-203, DOI: 10.5486/PMD.2014.5832.
%
\bibitem{M-Sz}
{Moln{\'a}r,~E.~--~Szirmai,~J.}
On homogeneous 3-geometries, balls and their optimal arrangements, especially 
in $\NIL$ and $\SOL$ spaces.
\emph{G-Slovak Journal for Geometry and Graphics} {\bf 19}  (2022), 37, 5-32.
%
\bibitem{stachel}
Moln{\'a}r,~E.~---~Szirmai,~J.
Packings with geodesic and translation balls and their visualizations in $\SLR$ space.
{\it Journal for Geometry and Graphics}  {\bf 26}(1) (2022), 51--64.
%
\bibitem{MSzV}
 {Moln{\'a}r,~E.~--~Szirmai,~J.~--~Vesnin,~A.}
 Projective metric realizations of cone-manifolds with singularities along 2-bridge knots and links.
 \emph{J. Geometry,}
{\bf 95} (2009), 91-133.
%
\bibitem{MSzV13}
 {Moln{\'a}r,~E.~--~Szirmai,~J.~--~Vesnin,~A.}
 Packings by translation balls in $\SLR$.
 \emph{J. Geometry,}
 {\bf 105 (2)} (2014), 287-306, DOI: 10.1007/s00022-013-0207-x.
%
\bibitem{R}
{Ratcliffe,~J.~G.}
\emph{Foundations of hyperbolic manifolds, (2nd ed.).}
{Graduate Texts in Mathematics 149. New York, NY: Springer.,} (2006).
 %
\bibitem{S}
{Scott,~P.}
The geometries of 3-manifolds.
\emph{Bull. London Math. Soc.}, {\bf15} (1983), 401--487.
  %
\bibitem{Sz07}
{Szirmai,~J.}
The densest geodesic ball packing by a type of $\NIL$ lattices.
\emph{Beitr. Algebra Geom.,}
 {\bf48(2)} (2007), 383--398.
 %
\bibitem{Sz10-2}
{Szirmai,~J.}
Geodesic ball packing in $\SXR$ space for generalized Coxeter space groups.
\emph{Beitr. Algebra Geom.,}
{\bf 52(2)} (2011), 413--430.
 %
\bibitem{Sz10-3}
{Szirmai,~J.}
Geodesic ball packing in $\HXR$ space for generalized Coxeter space groups.
\emph{Math. Commun.,}
{\bf 17/1} (2012), 151-170.
 %
\bibitem{Sz12-1}
{Szirmai,~J.}
Lattice-like translation ball packings in $\NIL$ space.
\emph{Publ. Math. Debrecen,}
{\bf 80/3-4} (2012), 427--440 DOI: 10.5486/PMD.2012.5117.
  %
\bibitem{Sz13-1}
{Szirmai,~J.}
Regular prism tilings in $\SLR$ space.
\emph{Aequat. Math.}, (2013), DOI 10.1007/s00010-013-0221-y.
%
\bibitem{Sz14}
Szirmai,~J
Non-periodic geodesic ball packings to infinite regular prism tilings in $\SLR$ space.  
\emph{Rocky Mountain J.  Math.}  {\bf 46/3} (2016), 1055--1070.
%
\bibitem{Sz13-2}
{Szirmai,~J.}
A candidate to the densest packing with equal balls in the Thurston geometries.
\emph{Beitr. Algebra Geom.,} (2013), DOI 10.1007/s13366-013-0158-2.
%
%
\bibitem{Sz202} Szirmai,~J.,
Interior angle sums of geodesic triangles in $\SXR$ and $\HXR$ geometries, 
{\it Bull. Academ. De Stiinte A Rep. Mol.}, {\bf{93}} Num 2 (2020), 44--61.
%
\bibitem{Sz22}
{Szirmai,~J.:}
Apollonius surfaces, circumscribed spheres of tetrahedra, Menelaus' and Ceva's theorems in $\SXR$ and $\HXR$ geometries,
\emph{Quarterly Journal of Mathematics}, {\bf 73} (2022), 477--494, doi: 10.1093/qmath/haab038, arXiv: 2012.06155.
%
\bibitem{Sz22-2}
{Szirmai,~J.:}
On Menelaus' and Ceva's theorems in $\NIL$ geometry,
\emph{Acta Univ. Sapientiae, Mathematica}, (to appear) (2023), arXiv: 2110.08877.
%
\bibitem{Sz22-3}
{Szirmai,~J.:}
Classical Notions and Problems in Thurston Geometries,
\emph{Submitted manuscript}, (2023), arXiv: 2203.05209.
%
\bibitem{T}
{Thurston,~W.~P.} (and {\sc Levy,~S.} editor)
\emph{Three-Dimensional Geometry and Topology.}
{Princeton University Press,} Princeton, New Jersey, Vol.{\bf 1} (1997).

 \end{thebibliography}
 \end{document}